\documentstyle[fullpage,amssymb,12pt]{article}
\title{Semi-Infinite Wedges and Vertex Operators}
\author{Eugene Stern}
\date{ \ }
\newtheorem{prop}{Proposition}[section]

\begin{document}
\maketitle
\begin{abstract}
The level 1 highest weight modules of the quantum affine algebra
$U_q(\widehat{\frak{sl}}_n)$ can be described as spaces of certain
semi-infinite wedges.  Using a $q$-antisymmetrization procedure,
these semi-infinite wedges can be realized inside an infinite
tensor product of evaluation modules.  This realization gives rise
to simple descriptions of vertex operators and (up to a scalar
function) their compositions.
\end{abstract}

\section{Representations of $\frak{sl}_{\infty}$}
\subsection{Infinite tensors and infinite wedges}
Let $V=\mbox{\bf C}^{\infty}$
be a countable-dimensional {\bf C}-vector space, with
basis $\{v_i\}_{i \in \mbox{\bf \scriptsize Z}}$.  The Lie algebra
$\frak{sl}_{\infty} = \frak{sl}_{\infty}(\mbox{\bf C})$,
consisting of infinite matrices with finitely
many non-zero entries and trace $0$, acts on $V$.  The elements
$e_i$, $f_i$, and $h_i$, where
\begin{eqnarray}
e_i \cdot v_j &=& \delta_{i,j-1} \cdot v_{j-1} \label{standard-1} \\
f_i \cdot v_j &=& \delta_{i,j} \cdot v_{j+1} \label{standard-2} \\
h_i \cdot v_j &=& (\delta_{i,j} - \delta_{i,j-1}) v_j, \label{standard-3}
\end{eqnarray}
generate $\frak{sl}_{\infty}$.

Consider the tensor product $V \otimes V \otimes V \otimes \cdots \ $.
In this tensor product, $\frak{sl}_{\infty}$ acts by
\begin{eqnarray}
e_i \cdot (v_{m_1} \otimes v_{m_2} \otimes \cdots ) &=& \sum_{j: \ m_j
= i+1} v_{m_1} \otimes \cdots \otimes v_{m_{j-1}} \otimes v_i \otimes
v_{m_{j+1}} \otimes \cdots \label{therm-1} \\
f_i \cdot (v_{m_1} \otimes v_{m_2} \otimes \cdots ) &=& \sum_{j: \ m_j
= i} v_{m_1} \otimes \cdots \otimes v_{m_{j-1}} \otimes v_{i+1} \otimes
v_{m_{j+1}} \otimes \cdots \\
h_i \cdot (v_{m_1} \otimes v_{m_2} \otimes \cdots ) &=& ( \# \{
r: m_r=i\} - \# \{r:m_r=i+1\}) \cdot
v_{m_1} \otimes v_{m_2} \otimes \cdots \label{therm-3}
\end{eqnarray}
To be precise, the infinite sums on the right hand side of
(\ref{therm-1})-(\ref{therm-3}) lie in an appropriate completion
of $V \otimes V \otimes V \otimes \cdots \ $.
Moreover, the action of $h_i$ given by (\ref{therm-3}) is defined
only for some tensors $v_{m_1} \otimes v_{m_2} \otimes \cdots \ $.
One restriction that certainly suffices is that each $v_i$ should
appear only finitely often among the $v_{m_j}$.

The infinite wedge product $\bigwedge^{\infty} V$ may be embedded
inside $V \otimes V \otimes V \otimes \cdots \ $ by an
antisymmetrization procedure.  Let $S_{\infty}$ denote
the infinite symmetric group, which is generated
by adjacent transpositions $\sigma_i = (i \ i+1)$, $i \in \mbox{\bf Z}^+$,
with the usual Coxeter relations.  ($S_{\infty}$ consists of bijections
$\mbox{\bf Z}^+ \to \mbox{\bf Z}^+$ which fix all but a finite number
of elements.)  The Bruhat length $l(\sigma)$ of an element $\sigma \in
S_{\infty}$ is the length $l$ of a minimal expression $\sigma =
\sigma_{i_1} \sigma_{i_2} \cdots \sigma_{i_l}$, where the $\sigma_{i_j}$
are adjacent transpositions.  (Thus, $(-1)^{l(\sigma)} = \mbox{sgn}
(\sigma)$.)

The antisymmetrization of $v_{m_1} \otimes v_{m_2} \otimes \cdots \ $ is
the pure wedge
\begin{equation}
v_{m_1} \wedge v_{m_2} \wedge \cdots = \sum_{\sigma \in S_{\infty}}
(v_{m_1} \otimes v_{m_2} \otimes \cdots ) \cdot (-1)^{l(\sigma)}
\sigma. \label{anti-sym}
\end{equation}
(Here $S_{\infty}$ acts on the right on infinite pure tensors in the
obvious way: $\sigma_i$ switches the $i$-th and $i+1$-st entries.)
Let $\bigwedge^{\infty} V \subseteq V \otimes V \otimes V \otimes
\cdots \ $ be the span of all pure wedges.

Since $v_{m_1} \wedge v_{m_2} \wedge \cdots \ $ is $0$ unless the $m_i$
are all distinct, the action of $\frak{sl}_{\infty}$ on the pure tensors
that appear in the expansion of any pure wedge is well defined.  Because
the actions of $\frak{sl}_{\infty}$ and $S_{\infty}$ commute, the
action of any $X \in \frak{sl}_{\infty}$ on the antisymmetrization of
$v_{m_1} \otimes v_{m_2} \otimes \cdots \ $ will yield the
antisymmetrization of $X \cdot (v_{m_1} \otimes v_{m_2} \otimes \cdots )$.
Thus the action of $\frak{sl}_{\infty}$ on tensors in $V \otimes V \otimes
V \otimes \cdots \ $ induces an action on $\bigwedge^{\infty} V$.

Consider the tensor
$$v_{(i)} = v_i \otimes v_{i-1} \otimes v_{i-2} \otimes \cdots \ .$$
Denote by $V_{(i)}$ the subspace of $V \otimes V \otimes V \otimes
\cdots \ $ spanned by pure tensors that are the same as $v_{(i)}$
after finitely many terms.  $V_{(i)}$ is preserved by the action of
both $\frak{sl}_{\infty}$ and $S_{\infty}$.

The pure wedge
$$v_{\Lambda_i} = v_i \wedge v_{i-1} \wedge v_{i-2} \wedge \cdots$$
lies in $V_{(i)}$.  It
is a highest weight vector of $\frak{sl}_{\infty}$
with highest weight $\Lambda_i$, and generates the irreducible
$\frak{sl}_{\infty}$-module $V_{\Lambda_i} \subseteq
V_{(i)}$ with highest weight $\Lambda_i$.  (Here $\Lambda_i$ is the
fundamental weight of $\frak{sl}_{\infty}$ defined by the equation
$\Lambda_i (h_j) = \delta_{ij}$.)
A basis for $V_{\Lambda_i}$ is given by wedges
$v_{m_1} \wedge v_{m_2} \wedge \cdots \ $ (with the $m_i$ decreasing),
which are the same as $v_{\Lambda_i}$ after finitely many
terms.  (See \cite{Kac} or \cite{KacRaina} for details.)  Such
wedges will be called {\em semi-infinite}.

\bigskip
\noindent
{\bf Remark} \hspace{1mm}
A semi-infinite wedge $v_{m_1} \wedge v_{m_2} \wedge \cdots
\subseteq V_{\Lambda_i}$ may be represented by a Young
diagram in the following way.  Set $\lambda_j
= m_j - (i-j+1)$, and notice that $m_j>m_{j+1}$ implies that $\lambda_j
\geq \lambda_{j+1}$.  Then set $\lambda=(\lambda_1, \lambda_2, \ldots)$.
After enough terms, the semi-infinite wedge becomes equal to
$v_{\Lambda_i}$,
so for large enough $j$, $m_j=i-j+1$, and $\lambda_j=0$.  In other words,
$\lambda$ is a finite Young diagram.  (For example, the highest weight
vector $v_{\Lambda_i}$ corresponds to the empty Young diagram.)
This correspondence sets up an isomorphism between representations
of $\frak{sl}_{\infty}$ on spaces of semi-infinite wedges and
representations on the space of Young diagrams written down in
\cite{paths}.

\subsection{Vertex operators}
\label{class-vertex}
If $v_{m_1} \wedge v_{m_2} \wedge \cdots \in V_{\Lambda_i}$,
then it is a sum of tensors all of which lie in $V_{(i)}$.  Given such
a tensor $v_{k_1} \otimes v_{k_2} \otimes v_{k_3} \otimes \cdots \in
V_{(i)}$, notice that the tensor $v_{k_2} \otimes v_{k_3} \otimes \cdots$
is an element of $V_{(i-1)}$.  This gives rise to a natural map
$\Phi_{(i)} : V_{(i)} \to V \otimes V_{(i-1)}$.  Since this map is
essentially the identity map, it commutes with the action of
$\frak{sl}_{\infty}$.

Now take the wedge $v_{m_1} \wedge v_{m_2} \wedge \cdots \ $
(remember, the $m_i$ are assumed to be decreasing)
and expand it as a sum of tensors.  Collect together the tensors
having $v_{m_j}$ as their first term.  The result is
$$(-1)^{j-1} v_{m_j} \otimes (v_{m_1} \wedge \cdots \wedge
v_{m_{j-1}} \wedge v_{m_{j+1}} \wedge \cdots ).$$
This shows that $\Phi_{(i)}$ maps $V_{\Lambda_i}$ into $V
\widehat{\otimes} V_{\Lambda_{i-1}}$.  (It is necessary to take an
appropriate completion since the image of a wedge in $V_{\Lambda_i}$
will be an infinite sum of $v_{m_j}$'s tensored with wedges in
$V_{\Lambda_{i-1}}$.)  In particular,
\begin{equation}
\Phi_{(i)} (v_{\Lambda_i}) = v_i \otimes v_{\Lambda_{i-1}} +
\sum_{j=1}^{\infty} (-1)^j v_{i-j} \otimes (v_i \wedge \cdots
\wedge v_{i-(j-1)} \wedge v_{i-(j+1)} \wedge \cdots ),
\end{equation}
i.e., $v_i$ is the ``matrix coefficient'' corresponding to
$v_{\Lambda_{i-1}}$.

Next, consider a composition
$$\Phi_{(i-(j-1))} \Phi_{(i-(j-2))} \cdots \Phi_{(i)} : V_{\Lambda_i}
\to \underbrace{V \otimes \cdots \otimes V}_{\mbox{\scriptsize $j$ times}}
\otimes V_{\Lambda_{i-j}}.$$
\begin{prop}
The matrix coefficient corresponding to $v_{\Lambda_{i-j}}$ is
$v_i \wedge v_{i-1} \wedge \cdots \wedge v_{i-(j-1)}$.
\end{prop}
{\em Proof.} \hspace{2mm}
Collect together terms ending in $v_{k_1} \otimes v_{k_2} \otimes
v_{k_3} \otimes \cdots \ $, where $k_1, k_2, k_3, \ldots$ is a
particular finite rearrangement of $i-j, i-j-1, i-j-2, \ldots \ $.

\section{Representations of $\widehat{\frak{sl}}_n$}
\label{affine-class}
\subsection{Evaluation modules}
The affine algebra $\widehat{\frak{sl}}_n$ has a standard {\em evaluation
representation} defined in the following way.  Let $E_i$, $F_i$, and $H_i$,
$i=0,1,\ldots ,n-1$, be the standard Serre generators of
$\widehat{\frak{sl}}_n$.  Consider an $n$-dimensional vector space
with basis $\{v_1, \ldots ,v_n\}$ on which these
generators act as follows:
\begin{eqnarray}
E_i \cdot v_j &=& \delta_{i,j-1} \cdot z^{\delta_{i,0}}
\cdot v_{j-1} \label{ev-mod-1} \\
F_i \cdot v_j &=& \delta_{i,j} \cdot z^{-\delta_{i,0}}
\cdot v_{j+1} \label{ev-mod-2} \\
H_i \cdot v_j &=& (\delta_{i,j} - \delta_{i+1,j}) \cdot
v_j \label{ev-mod-3}.
\end{eqnarray}
The indices in expressions (\ref{ev-mod-1})-(\ref{ev-mod-3}) should all
be read modulo $n$.  (For instance, if $j=1$, then $v_{j-1}=v_n$.)
It is easiest to regard $z$ as a formal variable by
tensoring over {\bf C} with the ring $\mbox{\bf C}[z,z^{-1}]$.  The
resulting $\widehat{\frak{sl}}_n$-module is denoted by $V(z)$.

$V(z)$ is related to the $\frak{sl}_{\infty}$-module $V=\mbox{\bf C}^{
\infty}$ of the
previous section in the following way.  Identify the basis $\{v_i\}_{i
\in \mbox{\bf \scriptsize Z}}$ of $V$ with the basis $\{z^j
\cdot v_i \}_{i \in 1, \ldots n, \ j \in \mbox{\bf \scriptsize Z}}$ of
$V(z)$ by $z^j \cdot v_i = v_{i-nj}$.  When $V(z)$ is
identified with $V$ in this way,
the generators $E_i$, $F_i$, and $H_i$ of $\widehat{\frak{sl}}_n$
act as infinite sums of the generators of $\frak{sl}_{\infty}$:
\begin{equation}
E_i = \sum_{j \equiv i \bmod n} e_j \hspace{.5in}
F_i = \sum_{j \equiv i \bmod n} f_j \hspace{.5in}
H_i = \sum_{j \equiv i \bmod n} h_j \label{correspond}
\end{equation}

\subsection{The thermodynamic limit}
As in the previous section, to build highest weight modules for
$\widehat{\frak{sl}}_n$, it is necessary to consider the infinite
tensor product
$$V_{z_1,z_2,z_3,\ldots} = V(z_1) \otimes V(z_2) \otimes V(z_3)
\otimes \cdots \ .$$
Unfortunately, it is not possible to define an honest action of
$\widehat{\frak{sl}}_n$ in this tensor product.  The most that can
be hoped for is a formal action of the Serre
generators $E_i$, $F_i$, and $H_i$, with the understanding that when
$E_i$ and $F_i$ act on a pure tensor, the result can be an {\em infinite}
sum of pure tensors.  Explicitly, let $E_i$ and $F_i$ act as formal sums
of operators
\begin{eqnarray}
\Delta(E_i) &=& \sum_{j=1}^{\infty} 1 \otimes \cdots \otimes 1
\otimes \underbrace{E_i}_{\mbox{\scriptsize $j$-th entry}} \otimes 1
\otimes 1 \otimes \cdots \\
\Delta(F_i) &=& \sum_{j=1}^{\infty} 1 \otimes \cdots \otimes 1
\otimes \underbrace{F_i}_{\mbox{\scriptsize $j$-th entry}} \otimes 1
\otimes 1 \otimes \cdots,
\end{eqnarray}
where the action of $E_i$ and $F_i$ in $V(z_j)$ is given by
(\ref{ev-mod-1})-(\ref{ev-mod-2}).
A highest or lowest weight vector is then just a vector killed
by all the $E_i$ or all the $F_i$.

An analogous formula for the action of $H_i$ will not work: if
defined naively, the action of $H_i$ on most tensors will
give a divergent answer. One way to get around this is to restrict
to a particular class of tensors
and to use a version of (\ref{correspond}).  Namely, $H_i$ acting on
$V(z)$ can be thought of as the sum over $d \in \mbox{\bf Z}$ of operators
$H_i(d)$, where
\begin{equation}
H_i(d) \cdot (z^{d'} v_j) = \delta_{d,d'} (\delta_{i,j} - \delta_{i+1,j})
\cdot z^{d'} v_j .
\end{equation}
Again, $v_0$ means $v_n$ here.
Then the action of $H_i(d)$ in $V_{z_1,z_2,z_3,\ldots}$ is given by the
infinite sum
\begin{equation}
\sum_{j=1}^{\infty} 1 \otimes \cdots \otimes 1
\otimes \underbrace{H_i(d)}_{\mbox{\scriptsize $j$-th entry}} \otimes 1
\otimes 1 \otimes \cdots . \label{class-coprod}
\end{equation}
Now, letting $H_i$ act on an arbitrary tensor in $V_{z_1,z_2,z_3\ldots}$ as
\begin{equation}
H_i = \sum_{d \in \mbox{\bf \scriptsize Z}} H_i(d) \label{mat-coeff}
\end{equation}
is still likely to give a divergent answer.  However, this formula can
at least be used for those tensors which eventually become periodic;
i.e., they have at their tail end an
infinite sequence of the form
$$z_k^d v_n \otimes z_{k+1}^d v_{n-1} \otimes \cdots \otimes
z_{k+n-1}^d v_1 \otimes z_{k+n}^{d+1} v_n \otimes z_{k+n+1}^{d+1} v_{n-1}
\otimes \cdots .$$
(All but finitely many of the $H_i(d)$ will act by $0$ on such
a tensor.)  Only tensors of this type (which will
be called {\em semi-infinite})
will be considered in the rest of this paper.

The action of the compositions $E_i F_i$ and $F_i E_i$ cannot be defined
in $V_{z_1,z_2,z_3\ldots}$: even the action on a semi-infinite tensor would
yield a divergent result.  However, the commutator of $E_i$ and $F_i$ can
still be taken, in the following formal way.  $E_i$ and $F_i$ can be
viewed as sums over $d \in \mbox{\bf Z}$ of operators $E_i(d)$ and
$F_i(d)$, whose action in $V(z)$ is given by
\begin{eqnarray}
E_i(d) \cdot (z^{d'} \cdot v_j) &=& \delta_{d,d'} \cdot \delta_{i,j-1}
\cdot z^{d'+\delta_{i,0}} \cdot v_{j-1} \label{E-comp} \\
F_i(d) \cdot (z^{d'} \cdot v_j) &=& \delta_{d,d'} \cdot \delta_{i,j}
\cdot z^{d'-\delta_{i,0}} \cdot v_{j+1}. \label{F-comp}
\end{eqnarray}
$E_i(d)$ and $F_i(d)$ act in $V_{z_1,z_2,z_3\ldots}$ as sums of the actions
in each component of the tensor product, in exact analogy with
(\ref{class-coprod}).  Then the equations
$$[E_i, F_i] = \left[ \sum_{d \in \mbox{\bf \scriptsize Z}} E_i(d) \ , \
\sum_{d \in \mbox{\bf \scriptsize Z}} F_i(d) \right] =
\sum_{d \in \mbox{\bf \scriptsize Z}} H_i(d) = H_i$$
hold formally in $V_{z_1,z_2,z_3\ldots}$ (more properly, they hold formally
in the subspace of $V_{z_1,z_2,z_3\ldots}$ spanned by semi-infinite tensors).

Notice that even though $V(z)$ and all finite tensor powers of it are
level 0 $\widehat{\frak{sl}}_n$-modules, the $H_i$ act in
$V_{z_1,z_2,z_3,\ldots}$ as though it were a level 1
$\widehat{\frak{sl}}_n$-module.
For example,
$$H_i \cdot (z_1 v_n \otimes \cdots \otimes z_n v_1 \otimes
z_{n+1}^2 v_n \otimes \cdots z_{2n}^2 v_1 \otimes \cdots) = \delta_{i,0}
\cdot (z_1 v_n \otimes \cdots \otimes z_n v_1 \otimes
z_{n+1}^2 v_n \otimes \cdots z_{2n}^2 v_1 \otimes \cdots),$$
because all $H_i(d)$ for $i \neq 0$, $d \neq 1$ act by $0$.  So
$z_1 v_n \otimes \cdots \otimes z_n v_1 \otimes
z_{n+1}^2 v_n \otimes \cdots z_{2n}^2 v_1 \otimes \cdots \ $
has weight $\Lambda_0$.  (In this section, $\Lambda_i$ is the $i$-th
fundamental weight of $\widehat{\frak{sl}}_n$, defined by $\Lambda_i
(H_j) = \delta_{i,j}$.)

It is worth pointing out that the formal action of the operators
$E_i$, $F_i$, $H_i \in \widehat{\frak{sl}}_n$
in the infinite tensor product is compatible with finite tensor
products.  That is, for $X = E_i, F_i$, or $H_i$,
\begin{eqnarray}
X \cdot (z_1^{d_1} v_{m_1} \otimes z_2^{d_2} v_{m_2} \otimes
z_3^{d_3} v_{m_3} \otimes \cdots ) &=& (X \cdot z_1^{d_1} v_{m_1})
\otimes (z_2^{d_2} v_{m_2} \otimes z_3^{d_3} v_{m_3} \otimes \cdots )
\label{compat} \\
& & + \ z_1^{d_1} v_{m_1} \otimes X \cdot (z_2^{d_2} v_{m_2} \otimes
z_3^{d_3} v_{m_3} \otimes \cdots) . \nonumber
\end{eqnarray}
(In (\ref{compat}), the action of $X$ on the left hand side and on the
second term on the right hand side are as defined above.)
This allows formal manipulation of the action in $V_{z_1,z_2,z_3,\ldots}$
as though it were an ordinary tensor product: it is legitimate to break
off finitely many factors, add on finitely many factors, etc.

As before, highest weight vectors are constructed by an antisymmetrization
procedure.  The antisymmetrization of $z_1^{d_1} v_{m_1} \otimes
z_2^{d_2} v_{m_2} \otimes \cdots \ $ is
\begin{equation}
z^{d_1} v_{m_1} \wedge z^{d_2} v_{m_2} \wedge \cdots = \sum_{\sigma
\in S_{\infty}} (z_1^{d_1} v_{m_1} \otimes z_2^{d_2} v_{m_2} \otimes
\cdots ) \cdot (-1)^{l(\sigma)} \sigma .
\end{equation}
Here $\sigma_i \in S_{\infty}$ acts by switching $v_{m_i}$ and $v_{
m_{i+1}}$ and by switching variables as well.  For example,
$$(z_1^{-5} v_3 \otimes z_2^3 v_6 \otimes z_3^{-2} v_1 \otimes \cdots)
\cdot \sigma_1 = z_1^3 v_6 \otimes z_2^{-5} v_3 \otimes z_3^{-2} v_1
\otimes \cdots \ .$$
The subscripts on the variables have been dropped in the infinite
wedge notation, but it should be understood that the antisymmetrization
is a sum of terms with $z_1$'s in the first factor, $z_2$'s in the
second factor, and so on.

Consider the wedge
\begin{equation}
v_{\Lambda_i} = zv_i \wedge zv_{i-1} \wedge \cdots \wedge
zv_1 \wedge z^2v_n \wedge \cdots \wedge z^2v_1 \wedge \cdots .
\label{class-hi-wt}
\end{equation}
Notice that this wedge is a sum of semi-infinite tensors.
Let $F_{(i)}$ denote the space spanned by wedges that are the same
as $v_{\Lambda_i}$ after finitely many terms.  As in the previous
section, such wedges will be called {\em semi-infinite}; they are
all sums of semi-infinite tensors.

The formal action of $E_i$, $F_i$, and $H_i \in
\widehat{\frak{sl}}_n$ on the semi-infinite tensors
in $V_{z_1,z_2,z_3,\ldots}$ defined above
induces an honest action on $F_{(i)}$.  The wedge $v_{\Lambda_i}$
is a highest weight vector of weight $\Lambda_i$: it is killed
by each $E_j$ because it is killed by all the $E_j(d)$.  (Since the
fundamental weights $\Lambda_i$ are usually indexed by $i=0,1,
\ldots ,n-1$, while the $v_i$ are indexed by $i=1,2,\ldots ,n$, it
is worth adding that $i=0$ on the left hand side of
(\ref{class-hi-wt}) corresponds to $i=n$ on the right hand side.)

The highest weight vector
$v_{\Lambda_i}$ generates an irreducible $\widehat{\frak{sl}}_n$-module
$V_{\Lambda_i} \subset F_{(i)}$ of
highest weight $\Lambda_i$.  (For a proof of irreducibity, see
\cite{KacRaina}.)  It is important to point out that $V_{\Lambda_i}$
is strictly smaller than $F_{(i)}$.  For example, for $n=2$ and
$i=0$, the wedges
$v_2 \wedge z v_1 \wedge z^2 v_2 \wedge z^2 v_1 \wedge \cdots$ and
$v_1 \wedge z v_2 \wedge z^2 v_2 \wedge z^2 v_1 \wedge \cdots \in
F_{(0)}$ do not lie in $V_{\Lambda_0}$, although their sum does, being
$F_1 F_0 \cdot v_{\Lambda_0}$.  On the other hand, $F_{(0)}$ is a
unitary $\widehat{\frak{sl}}_n$-module (see \cite{KacRaina}),
and therefore completely reducible, so there is a projection $p_i: F_{(i)}
\to V_{\Lambda_i}$ which kills all the components of $F_{(i)}$ except
the one generated by $v_{\Lambda_i}$.  In particular, $p_i(v_{\Lambda_i})
=v_{\Lambda_i}$.

Once again, there is a correspondence between wedges
and Young diagrams.  For example, if $i=0$, the semi-infinite wedge
$$  v_3 \wedge v_1
\wedge zv_{n-2} \wedge zv_{n-3} \wedge \ldots \wedge zv_1 \wedge
z^2v_n \wedge \cdots \wedge z^2v_1 \wedge \cdots $$
corresponds to the Young diagram $(3,2)$.  In the notation of
\cite{paths}, the value of $i$ (i.e., the highest weight)
determines the way in which the diagram is to be colored.

\subsection{Vertex operators}
\label{vertex}
Just as in the first section, splitting off the first component
of a tensor defines intertwiners; in this case,
$$\widetilde{\Phi}_{(i)} : F_{(i)} \to V(z) \ \widehat{\otimes}
\ F_{(i-1)}.$$
It is worthwhile to say explicitly what this means in the cases
$i=0$ and $i=1$.  $\widetilde{\Phi}_{(0)}$
maps $F_{(0)}$ into $V(z) \ \widehat{\otimes} \ F_{(n-1)}$.
(I.e., the indices should really be read modulo $n$.)
$\widetilde{\Phi}_{(1)}$ maps $F_{(1)}$, spanned by wedges that
look like $zv_1 \wedge z^2v_n \wedge \cdots
\wedge z^2v_1 \wedge \cdots \ $ after finitely many terms, to $V(z) \
\widehat{\otimes} \ W$, where $W$ is spanned by wedges that look like
$z^2v_n \wedge z^2v_{n-1} \wedge \cdots
\wedge z^2v_1 \wedge z^3v_n \wedge \cdots \ $ after finitely many
terms.  $W$ is evidently
isomorphic to, and can be identified with, $F_{(0)}$. (That is,
the $\widetilde{\Phi}_{(i)}$ cycle around after $n$ iterations:
$\widetilde{\Phi}_{(1)}$ is followed by $\widetilde{\Phi}_{(0)}$,
which acts on the $\widehat{\frak{sl}}_n$-module spanned by
wedges that look like $z^2 v_n \wedge z^2 v_{n-1} \wedge \cdots
\wedge z^2 v_1 \wedge z^3 v_n \wedge \cdots $ after finitely many
terms, and so on.)  Each $\widetilde{\Phi}_{(i)}$ is an intertwiner
because of (\ref{compat}).

\begin{prop}
When the composition
$$\widetilde{\Phi}_{(n-(j-1))} \widetilde{\Phi}_{(n-(j-2))} \cdots
\widetilde{\Phi}_{(0)} (\widetilde{\Phi}_{(1)} \widetilde{\Phi}_{(2)}
\cdots \widetilde{\Phi}_{(0)} )^d :
F_{(0)} \to V(z_1) \otimes V(z_2) \otimes \cdots \otimes V(z_{nd+j})
\ \widehat{\otimes} \  F_{(n-j)}$$
acts on the highest weight vector $v_{\Lambda_0}$,
the matrix coefficient corresponding to $v_{\Lambda_{n-j}}$ is
\begin{equation}
z v_n \wedge \cdots \wedge z v_1 \wedge \cdots \wedge z^d v_n \wedge
\cdots \wedge z^d v_1 \wedge z^{d+1} v_n \wedge \cdots \wedge z^{d+1}
v_{n-(j-1)}. \label{high-to-high}
\end{equation}
\end{prop}
The proof is as in the previous section.

Next, define the intertwiner $\Phi_{(i)} = (\mbox{\rm id} \otimes p_{i-1})
\circ \widetilde{\Phi}_{(i)} \circ j_i : V_{\Lambda_i} \to V(z) \
\widehat{\otimes} \ V_{\Lambda_{i-1}}$.  Here $j_i: V_{\Lambda_i} \to
F_{(i)}$ is the inclusion map, and $p_{i-1}: F_{(i-1)} \to
V_{\Lambda_{i-1}}$
is the projection.  The intertwiners $\Phi_{(i)}$ are examples of vertex
operators (see, for example, \cite{qKZ}).

Consider the composition of vertex operators
\begin{equation}
\Phi_{(n-(j-1))}
\Phi_{(n-(j-2))} \cdots \Phi_{(0)} (\Phi_{(1)} \Phi_{(2)} \cdots
\Phi_{(0)} )^d. \label{composition-1}
\end{equation}
This is an intertwiner from $V_{\Lambda_0}$ to
$V(z_1) \otimes \cdots \otimes V(z_{nd + j}) \ \widehat{\otimes}
\ V_{\Lambda_{n-j}}$,
whose action can be computed up to scalars as follows.  Notice that the
composition
\begin{equation}
p_{n-j} \ \circ \ \widetilde{\Phi}_{(n-(j-1))}
\widetilde{\Phi}_{(n-(j-2))} \cdots
\widetilde{\Phi}_{(0)} (\widetilde{\Phi}_{(1)} \widetilde{\Phi}_{(2)}
\cdots \widetilde{\Phi}_{(0)} )^d \ \circ \ j_0 \label{composition-2}
\end{equation}
is also an intertwiner
from $V_{\Lambda_0}$ to $V(z_1) \otimes \cdots \otimes V(z_{nd + j})
\ \widehat{\otimes} \ V_{\Lambda_{n-j}}$,
and since the space of such intertwiners is one-dimensional, the two
are equal up to scalars.  (Neither one is zero:
(\ref{composition-1}) is a non-zero composition of vertex operators,
and (\ref{composition-2}) is non-zero since its action on $v_{\Lambda_0}$
is non-zero.)  This scalar depends on the ratios $z_k/z_l$.
In other words, to compute an iteration of vertex operators up to
scalars, it suffices to iterate the $\widetilde{\Phi}_{(i)}$
and take one projection at the end, instead of applying the projections
at every step.  Together with the preceding proposition, this implies

\begin{prop}
When the composition (\ref{composition-1}) acts on the highest weight
vector $v_{\Lambda_0}$,
the matrix coefficient corresponding to $v_{\Lambda_{n-j}}$ is
a scalar multiple of (\ref{high-to-high}). \label{vertex-comp}
\end{prop}
It is well-known that the ``highest-to-highest matrix coefficients'' of
a composition of vertex operators are given by the Knizhnik-Zamolodchikov
equations (see, for example, \cite{qKZ}).  Thus, appropriate multiples
of the wedges given by (\ref{high-to-high}) should be solutions
to these equations.  As an example of this, consider
an iteration $(\Phi_{(n-1)} \cdots \Phi_{(1)} \Phi_{(0)})^N$.
Then the KZ system is a system of $nN$ differential equations with
values in the $\frak{sl}_n$-module $V^{\otimes nN}$ ($V=\mbox{\bf C}^n$):
\begin{equation}
\left( \frac{\partial}{\partial z_i} - \frac1{n+1} \sum_{j \neq i}
\frac{t_{ij}}{z_i - z_j} \right) \cdot F = 0.  \label{KZeq}
\end{equation}
Here $t_{ij} = P_{ij} - \ \frac1n \ \mbox{\rm id}$, where $P_{ij}$ acts on
$V^{\otimes nN}$ by switching the vectors (but not the variables) in
the $i$-th and $j$-th places.  Then the solution is indeed a multiple of
(\ref{high-to-high}).  It is given by the following proposition, which
came out of discussions with Nicolai Reshetikhin:

\begin{prop}
The function
\begin{eqnarray*}
F(z_1, \ldots , z_{nN}) &=& \prod_{1 \leq i<j \leq nN} (z_i - z_j)^{-1/n}
\cdot \\
& & \hspace{.8in}
(v_n \wedge v_{n-1} \wedge \cdots \wedge v_1 \wedge zv_n \wedge \cdots
\wedge z^{N-1} v_n \wedge \cdots \wedge z^{N-1} v_1) \\
\smallskip
 &=& \prod_{1 \leq i<j \leq nN} (z_i - z_j)^{-1/n} z_1^{-1} z_2^{-1}
\cdots z_{nN}^{-1} \cdot \\
& & \hspace{.8in}
(z v_n \wedge z v_{n-1} \wedge \cdots \wedge z v_1 \wedge z^2v_n \wedge
\cdots \wedge z^N v_n \wedge \cdots \wedge z^N v_1)
\end{eqnarray*}
is a solution to the $\frak{sl}_n$-KZ system.
\end{prop}

\bigskip
\noindent
{\bf Remarks}

\smallskip
\noindent
1.  The highest weight vector (i.e., the infinite antisymmetrization)
$v_{\Lambda_0} \in V_{z_1,z_2,z_3\ldots}$
should be regarded essentially as a limit as $N \to \infty$ of the
function $F(z_1,\ldots , z_{nN})$ (the finite antisymmetrization)
given above.

\smallskip
\noindent
2.  All iterations of vertex operators in this section began with
$V_{\Lambda_0}$ purely for convenience.  The results have obvious analogs
for iterations of vertex operators beginning with any level 1 highest
weight $\widehat{\frak{sl}}_n$-module.

\section{Representations of $U_q(\frak{sl}_{\infty})$}
\label{sl-inf}
\subsection{Preliminaries}
The main subject of this paper is a quantization of the picture laid
out in the first two sections.  To begin with,
there is a quantization of $\frak{sl}_{\infty}$, denoted by
$U_q(\frak{sl}_{\infty})$, which also acts on $V=\mbox{\bf C}^{\infty}$.
$U_q(\frak{sl}_{\infty})$ is generated by elements $e_i$, $f_i$,
$k_i$, and $k_i^{-1}$, $i \in \mbox{\bf Z}$, with relations
\begin{eqnarray}
k_i k_j &=& k_j k_i \label{begingen} \\
k_i e_i &=& q e_i k_i \\
k_i f_i &=& q^{-1} f_i k_i \\
e_i f_j - f_j e_i &=& \delta_{i,j} \frac{k_i - k_i^{-1}}
  {q - q^{-1}} \label{commutator} \\
e_i e_j &=& e_j e_i \hspace{.8in} \mbox{\rm if}
          \hspace{2mm} |i-j|>1 \\
f_i f_j &=& f_j f_i \hspace{.8in} \mbox{\rm if}
          \hspace{2mm} |i-j|>1 \\
0 &=& e_i^2 e_{i \pm 1} - (q + q^{-1})e_i e_{i \pm 1} e_i
        + e_{i \pm 1} e_i^2 \\
0 &=& f_i^2 f_{i \pm 1} - (q + q^{-1})f_i f_{i \pm 1} f_i
        + f_{i \pm 1} f_i^2. \label{endgen}
\end{eqnarray}
$U_q(\frak{sl}_{\infty})$ acts on $V$ as follows: the
action of $e_i$ and $f_i$ is the same as in the classical case
(see (\ref{standard-1})-(\ref{standard-2})), while
$k_i$ acts as $q^{h_i}$, where the action of $h_i$ is given by
(\ref{standard-3}).

There is a coproduct on $U_q(\frak{sl}_{\infty})$ given by
\begin{eqnarray}
\Delta(k_i) &=& k_i \otimes k_i \\
\Delta(e_i) &=& e_i \otimes k_i + 1 \otimes e_i \\
\Delta(f_i) &=& f_i \otimes 1 + k_i^{-1} \otimes f_i .
\end{eqnarray}
This coproduct gives rise to the following action of $U_q(\frak{sl}_{
\infty})$ on certain infinite pure tensors $v_{m_1} \otimes v_{m_2}
\otimes \cdots \ $:
\begin{eqnarray}
e_i \cdot (v_{m_1} \otimes v_{m_2} \otimes \cdots) &=& \sum_{j: \, m_j
= i+1} q^{\# \{r:\, r>j, \, m_r=i\} - \# \{r: \, r>j, \, m_r=i+1 \}}_{
\textstyle \hspace{.4in} \cdot \hspace{.1in}
v_{m_1} \otimes \cdots \otimes v_{m_{j-1}} \otimes
v_i \otimes v_{m_{j+1}} \otimes \cdots} \label{coprodaction1} \\
f_i \cdot (v_{m_1} \otimes v_{m_2} \otimes \cdots) &=& \sum_{j: \, m_j=i}
q^{\# \{r:\, r<j, \, m_r=i+1\} - \# \{r: \, r<j, \, m_r=i \}}_{
\textstyle \hspace{.4in} \cdot \hspace{.1in}
v_{m_1} \otimes \cdots \otimes v_{m_{j-1}} \otimes
v_{i+1} \otimes v_{m_{j+1}} \otimes \cdots} \label{coprodaction2} \\
k_i \cdot (v_{m_1} \otimes v_{m_2} \otimes \cdots ) &=& q^{
\# \{r: \, m_r=i \} - \# \{r: \, m_r=i+1 \}} \cdot
v_{m_1} \otimes v_{m_2} \otimes \cdots \ . \label{coprodaction3}
\end{eqnarray}
As before, the terms on the right hand side
lie in an appropriate completion of the infinite tensor product.
To ensure the action is well defined, the domain is again restricted
to the subspace of $V \otimes V \otimes V \otimes \cdots \ $ spanned by
tensors $v_{m_1} \otimes v_{m_2} \otimes \cdots \ $ in which all $v_j$
appear only finitely many times.

\subsection{$q$-antisymmetrization}
For $i \in \mbox{\bf Z}$, consider the pure tensor
$$v_{(i)} = v_i \otimes v_{i-1} \otimes v_{i-2} \otimes \cdots.$$
Denote by $V_{(i)}$ the subspace of $V \otimes V \otimes V \otimes
\cdots$ spanned by all pure tensors that
are the same as $v_{(i)}$ after finitely many terms.  Notice that
the action of $U_q(\frak{sl}_{\infty})$ on each $V_{(i)}$ is well
defined.  $v_{(i)}$ has weight $\Lambda_i$ with respect to the
subalgebra $U_q(\frak{h}) \subseteq
U_q(\frak{sl}_{\infty})$ generated by $\{k_i, k_i^{-1}
\}_{i \in \mbox{\bf \scriptsize Z}}$.

Consider the $q$-antisymmetrization of $v_{(i)}$ given by
\begin{equation}
v_{\Lambda_i} = \sum_{\sigma \in S_{\infty}} v_{(i)} \cdot
(-q)^{l(\sigma)} \sigma. \label{high-wt}
\end{equation}
$v_{\Lambda_i}$ is an infinite sum of elements of $V_{(i)}$, all
of which have weight $\Lambda_i$.
\begin{prop}
Every $e_j$ acts on $v_{\Lambda_i}$ by zero.
\end{prop}
{\em Proof.} \hspace{2mm}
For concreteness, assume that $i=0$.  In this case, if $j \geq
0$, there is nothing to prove.  If $j<0$, partition $S_{\infty}$
into left cosets of the subgroup $H=\{ \mbox{id}, \sigma_{-j} \}$.  Each
such coset looks like $\{\sigma, \sigma_{-j} \sigma\}$, where
$l(\sigma_{-j} \sigma) = l(\sigma)+1$.  Group together terms in
(\ref{high-wt}) in pairs corresponding to these cosets, and
consider one such pair.  The term
corresponding to $\sigma$ has $v_{j+1}$ appearing to the left of $v_j$.
The term
corresponding to $\sigma_{-j} \sigma$ looks just like the term
corresponding to $\sigma$ except that it has the opposite sign,
an extra factor of $q$, and $v_{j+1}$ and $v_j$ are switched, so that
$v_{j+1}$ appears to the right of $v_j$.  By (\ref{coprodaction1}),
$e_j$ kills the sum of these two terms.  All other pairs are killed in
the same way.

\medskip
\noindent
If $v_{m_1} \otimes v_{m_2} \otimes \cdots \in V_{(i)}$ is such that
the $m_i$ are decreasing, set
\begin{equation}
v_{m_1} \wedge_q v_{m_2} \wedge_q \cdots = \sum_{\sigma \in S_{\infty}}
(v_{m_1} \otimes v_{m_2} \otimes \cdots ) \cdot (-q)^{l(\sigma)} \sigma .
\label{q-antisym}
\end{equation}
The $q$-antisymmetrized tensor
$v_{m_1} \wedge_q v_{m_2} \wedge_q \cdots \ $ will be called a
{\em $q$-wedge}.
Let $V_{\Lambda_i}$ denote the space spanned by the $q$-wedges
$v_{m_1} \wedge_q v_{m_2} \wedge_q \cdots \ $ (with the $m_i$
decreasing), which are the same as $v_{\Lambda_i}$ after finitely
many terms.  (Such $q$-wedges will be called {\em semi-infinite}.)

By the above, $v_{\Lambda_i}$ generates a highest weight $U_q(\frak{sl}_{
\infty})$-submodule in $V_{(i)}$, of highest weight $\Lambda_i$.  This
module is spanned by the semi-infinite $q$-wedges described above, and
$U_q(\frak{sl}_{\infty})$ acts on these wedges in the obvious way.
(Notice that $f_j$ kills any $q$-wedge of the form $ \ \cdots \wedge_q
v_{j+1} \wedge_q v_j \wedge_q \cdots \ $, which means that the action of
$U_q(\frak{sl}_{\infty})$ cannot generate any $q$-wedges in which any
$v_j$ appears more than once.)  This
representation is irreducible because it is a $q$-deformation of the
irreducible $\frak{sl}_{\infty}$-module $V_{\Lambda_i}$ constructed in
the first section.

Identifying $q$-wedges with Young diagrams as before gives rise to an
isomorphism of $V_{\Lambda_i}$ with the level one highest weight
$U_q(\frak{sl}_{\infty})$-modules described in \cite{crystal}.

\subsection{Vertex operators}
Intertwiners $\Phi_{(i)}: V_{\Lambda_i} \to V \widehat{\otimes}
V_{\Lambda_{i-1}}$ can be defined exactly as in the classical case,
by splitting off the first component of every tensor.  The $q$-analogs
of the results of Section \ref{class-vertex} are the following:
\begin{prop}
The image of the highest weight vector $v_{\Lambda_i} \in V_{\Lambda_i}$
under $\Phi_{(i)}$ is given by
\begin{equation}
\Phi_{(i)} (v_{\Lambda_i}) = v_i \otimes v_{\Lambda_{i-1}} +
\sum_{j=1}^{\infty} (-q)^j v_{i-j} \otimes (v_i \wedge_q \cdots
\wedge_q v_{i-(j-1)} \wedge_q v_{i-(j+1)} \wedge_q \cdots ).
\end{equation}
\end{prop}
\begin{prop}
Under the composition
$$\Phi_{(i-(j-1))} \Phi_{(i-(j-2))} \cdots \Phi_{(i)} : V_{\Lambda_i}
\to \underbrace{V \otimes \cdots \otimes V}_{\mbox{\scriptsize $j$
{\em times}}} \otimes V_{\Lambda_{i-j}},$$
the matrix coefficient corresponding to
$v_{\Lambda_{i-j}}$ is $v_i \wedge_q v_{i-1} \wedge_q \cdots \wedge_q
v_{i-(j-1)}$.
\end{prop}

\section{Representations of $U_q(\widehat{\frak{sl}}_n)$}
\subsection{Evaluation modules}
The quantum affine algebra $U_q(\widehat{\frak{sl}}_n)$ is an algebra
generated by elements $E_i$, $F_i$, $K_i$, and $K_i^{-1}$ for $i=0,1,
\ldots , n-1$.  These elements satisfy the analogs of relations
(\ref{begingen})-(\ref{endgen}), for example,
\begin{equation}
E_i F_j - F_j E_i = \delta_{i,j} \frac{K_i - K_i^{-1}}
  {q - q^{-1}}, \label{q-aff-comm}
\end{equation}
with the added stipulation that the indices in all the relations
should be read modulo $n$.

As in the classical case, $U_q(\widehat{\frak{sl}}_n)$ acts on an
evaluation module $V(z)$.  Precisely, if $V(z)$ is the same vector
space as in Section \ref{affine-class}, then the generators of
$U_q(\widehat{\frak{sl}}_n)$ act as follows:
\begin{eqnarray}
K_i \cdot v_j &=& q^{(\delta_{i,j} - \delta_{i+1,j})} \cdot v_j
\label{eval-mod-1} \\
E_i \cdot v_j &=& \delta_{i,j-1} \cdot z^{\delta_{i,0}} \cdot v_{j-1} \\
F_i \cdot v_j &=& \delta_{i,j} \cdot z^{-\delta_{i,0}} \cdot v_{j+1}
\label{eval-mod-3} .
\end{eqnarray}
Here again, all the indices should be read modulo $n$.  So $E_i$ and
$F_i$ act as in Section \ref{affine-class}, while $K_i$ acts as
$q^{H_i}$.  Under the identification $z^j \cdot v_i = v_{i-nj}$ of
$V(z)$ with $\mbox{\bf C}^{\infty}$, $E_i$, $F_i$, and $K_i$ acting in
$V(z)$ can be expressed in terms of elements of $U_q(\frak{sl}_{\infty})$
acting in $\mbox{\bf C}^{\infty}$:
\begin{equation}
E_i = \sum_{j \equiv i \bmod n} e_j \hspace{.5in}
F_i = \sum_{j \equiv i \bmod n} f_j \hspace{.5in}
K_i = \prod_{j \equiv i \bmod n} k_j \label{q-correspond}
\end{equation}
The identification of these two modules should not be taken to mean that
$U_q(\widehat{\frak{sl}}_n)$ can in general be considered to be sitting
inside $U_q(\frak{sl}_{\infty})$.  The fact that the operators on
$\mbox{\bf C}^{\infty}$ defined by equations (\ref{q-correspond}) satisfy
the relation (\ref{q-aff-comm}) is a consequence of the following equation
for $U_q(\frak{sl}_{\infty})$ acting in $\mbox{\bf C}^{\infty}$:
$$\sum_{j \equiv i \bmod n} \frac{k_j - k_j^{-1}}{q-q^{-1}} =
\frac{1}{q-q^{-1}} \left( \prod_{j \equiv i \bmod n} k_j \ -
\prod_{j \equiv i \bmod n} k_j^{-1} \right).$$

\subsection{The thermodynamic limit}
To build highest weight modules, it is again necessary to consider the
infinite tensor product
$$V_{z_1, z_2, z_3, \ldots} = V(z_1) \otimes V(z_2) \otimes V(z_3)
\otimes \cdots \ .$$
As in the previous section, a coproduct is needed to define how the
generators of $U_q(\widehat{\frak{sl}}_n)$ act (at least formally)
on an appropriate subspace of $V_{z_1,z_2,z_3,\ldots}$.
The coproduct on $U_q(\widehat{\frak{sl}}_n) $ is analogous to the
one on $U_q(\frak{sl}_{\infty})$; explicitly it is given by
\begin{eqnarray}
\Delta(K_i) &=& K_i \otimes K_i  \label{fin-coprod-1} \\
\Delta(E_i) &=& E_i \otimes K_i + 1 \otimes E_i \label{fin-coprod-2} \\
\Delta(F_i) &=& F_i \otimes 1 + K_i^{-1} \otimes F_i . \label{fin-coprod-3}
\end{eqnarray}
Iterating this gives rise to the ``infinite coproduct''
\begin{eqnarray}
\Delta^{\infty} (K_i) &=& K_i \otimes K_i \otimes K_i \otimes \cdots
\label{coprod-1} \\
\Delta^{\infty} (E_i) &=& \sum_{j=1}^{\infty} 1 \otimes \cdots \otimes 1
\otimes \underbrace{E_i}_{\mbox{\scriptsize $j$-th entry}} \otimes K_i
\otimes K_i \otimes \cdots \label{coprod-2} \\
\Delta^{\infty} (F_i) &=& \sum_{j=1}^{\infty} K_i^{-1} \otimes \cdots
\otimes K_i^{-1}
\otimes \underbrace{F_i}_{\mbox{\scriptsize $j$-th entry}} \otimes 1
\otimes 1 \otimes \cdots \label{coprod-3}
\end{eqnarray}
This coproduct should define a formal action of the operators $K_i$,
$E_i$, and $F_i \in U_q(\widehat{\frak{sl}}_n)$ on semi-infinite tensors
(semi-infinite meaning the same thing here as in Section 2) in
$V_{z_1,z_2,z_3\ldots}$.  This is done exactly as in Section 2;
the only new feature of the quantum case is that it is necessary to say how
$K_i \otimes K_i \otimes K_i \otimes \cdots \ $ acts. The only possibility
is to make it act as $q^{H_i}$, where the action of $H_i$ in $V_{z_1,z_2,
z_3,\ldots}$ is given by (\ref{mat-coeff}).  For example,
\begin{eqnarray*}
\lefteqn{(K_i \otimes K_i \otimes K_i \otimes \cdots )
\cdot (z_1 v_n \otimes \cdots \otimes z_n v_1 \otimes
z_{n+1}^2 v_n \otimes \cdots z_{2n}^2 v_1 \otimes \cdots) =} \\
& & \hspace{2in} q^{\delta_{i,0}}
\cdot (z_1 v_n \otimes \cdots \otimes z_n v_1 \otimes
z_{n+1}^2 v_n \otimes \cdots \otimes z_{2n}^2 v_1 \otimes \cdots).
\end{eqnarray*}
With these definitions, the operators $K_i$, $E_i$, and $F_i$ act
formally on semi-infinite tensors in $V_{z_1,z_2,z_3,\ldots}$.  This
formal action is compatible with finite tensor
products (i.e., an analog of (\ref{compat}) holds).

\subsection{$q$-antisymmetrization}
Some care is required to produce highest weight vectors for
$U_q(\widehat{\frak{sl}}_n)$, since the naive $q$-antisymmetrization
given by equation (\ref{high-wt}) does not work. The correct approach
is to use a form of quantum Weyl duality. The symmetric group $S_d$
has a quantum analog known as a Hecke algebra, to be denoted here by
$H_d(q^2)$.  $H_d(q^2)$ is an $d!$-dimensional algebra generated by
elements $T_i$, $i=1,\ldots ,d-1$, satisfying the relations
\begin{eqnarray}
T_i^2 &=& (q^2-1) \ T_i + q^2 \label{heckegen-1} \\
T_i T_{i+1} T_i &=& T_{i+1} T_i T_{i+1} \label{heckegen-2} \\
T_i T_j &=& T_j T_i  \hspace{.8in} \mbox{\rm if} \hspace{2mm} |i-j|>1 .
\label{heckegen-3}
\end{eqnarray}
The elements $T_i$ are $q$-analogs of the adjacent transpositions
$\sigma_i = (i \ i+1)$ in the symmetric group $S_d$.

$H_d(q^2)$ acts on the right on the tensor product $V(z_1) \otimes \cdots
\otimes V(z_d)$ as follows.  Write elements $z_1^{j_1} v_{m_1}
\otimes \cdots
\otimes z_d^{j_d} v_{m_d}$ as $(v_{m_1} \otimes \cdots \otimes v_{m_d})
\cdot z_1^{j_1} \cdots z_d^{j_d}$.  $S_d$ can act on both the tensor part
and the polynomial part of such an expression.  The action on the tensor
part is the usual one, permuting factors:
$$(v_{m_1} \otimes \cdots \otimes v_{m_d})^{\sigma_i} =
v_{m_1} \otimes \cdots \otimes v_{m_{i+1}} \otimes v_{m_i} \otimes
\cdots \otimes v_{m_d}.$$
Similarly, the action on the polynomial part is to permute variables:
if $\mbox{\bf z} = z_1^{j_1} \cdots z_d^{j_d}$, then
$$\mbox{\bf z}^{\sigma_i} = (z_1^{j_1} \cdots z_d^{j_d})^{\sigma_i}
=z_1^{j_1} \cdots z_i^{j_{i+1}} z_{i+1}^{j_i} \cdots z_d^{j_d}.$$
Then,
\begin{equation}
((v_{m_1} \otimes \cdots \otimes v_{m_d}) \cdot \mbox{\bf z})
\cdot (T_i) =
\left\{
\begin{array}{ll}
- q(v_{m_1} \otimes \cdots \otimes v_{m_d})^{\sigma_i} \cdot
\mbox{\bf z}^{\sigma_i} \\
\hspace{3mm} - (q^2-1)(v_{m_1} \otimes \cdots \otimes v_{m_d}) \cdot
\frac{z_{i+1} \mbox{\bf \footnotesize z}^{\sigma_i} -
z_i \mbox{\bf \footnotesize z} }{z_i-z_{i+1}} & \hspace{3mm}
\mbox{if $m_i<m_{i+1}$} \\
- (v_{m_1} \otimes \cdots \otimes v_{m_d}) \cdot \mbox{\bf z}^{\sigma_i} \\
\hspace{3mm} - (q^2-1) (v_{m_1} \otimes \cdots \otimes v_{m_d}) \cdot
\frac{z_i(\mbox{\bf \footnotesize z}^{\sigma_i} -
\mbox{\bf \footnotesize z})}{z_i-z_{i+1}} & \hspace{3mm}
\mbox{if $m_i=m_{i+1}$} \\
- q(v_{m_1} \otimes \cdots \otimes v_{m_d})^{\sigma_i} \cdot
\mbox{\bf z}^{\sigma_i} \\
\hspace{3mm} - (q^2-1)(v_{m_1} \otimes \cdots \otimes v_{m_d}) \cdot
\frac{z_i(\mbox{\bf \footnotesize z}^{\sigma_i} -
\mbox{\bf \footnotesize z})}{z_i-z_{i+1}} & \hspace{3mm}
\mbox{if $m_i>m_{i+1}$}
\end{array} \right. \label{hecke-action}
\end{equation}

\medskip
\noindent
{\bf Remarks}

\smallskip
\noindent
1.  Notice that, for example, $(v_i \otimes v_i) \cdot T_1 = -v_i
\otimes v_i$;
also, $(v_{i+1} \otimes v_i) \cdot T_1 = -q v_i \otimes v_{i+1}$,
and $(v_1 \otimes z_2 v_n) \cdot T_1 = -q z_1 v_n \otimes v_1$.  (These
equations remain true if the left and right hand sides are multiplied by
$(z_1 z_2)^d$ for any $d$.)  In particular, the $T_i's$, rather than the
$-T_i$'s, will be used to do $q$-antisymmetrization.

\smallskip
\noindent
2.  The {\em affine} Hecke algebra $\widehat{H}_d(q^2)$ acts on $V(z_1)
\otimes
\cdots \otimes V(z_d)$ as the centralizer of the action of $U_q(\widehat{
\frak{sl}}_n)$ given by the coproduct in
(\ref{fin-coprod-1})-(\ref{fin-coprod-3}) (equivalently, the action of
$U_q(\widehat{\frak{sl}}_n)$ is
given by a finite version of (\ref{coprod-1})-(\ref{coprod-3})).  The
action of $H_d(q^2)$ written down above comes from regarding it as a
subalgebra of $\widehat{H}_d(q^2)$ in the obvious way.

\bigskip
There is a chain of inclusions $H_1(q^2) \subset H_2(q^2) \subset H_3(q^2)
\subset \cdots \ $, so equation (\ref{hecke-action}) also defines an
action of the {\em infinite} Hecke algebra $H_{\infty}(q^2)=\bigcup_{d
\geq 1} H_d(q^2)$ on the thermodynamic limit $V_{z_1,z_2,z_3,\ldots}$.
($H_{\infty}(q^2)$ is generated by elements $T_1, T_2, T_3,\ldots$ with
the corresponding relations.)  This action commutes with the action
of $U_q(\widehat{\frak{sl}}_n)$ in $V_{z_1,z_2,z_3,\ldots}$ since any
element acts in only finitely many factors.

There is a natural basis for $H_{\infty}(q^2)$ made up of elements
$T_{\sigma}$ corresponding to $\sigma \in S_{\infty}$.  More precisely,
if $\sigma = \sigma_{i_1} \sigma_{i_2} \cdots \sigma_{i_l}$ is a minimal
length expansion of $\sigma \in S_{\infty}$ in terms of adjacent
transpositions, let $T_{\sigma} = T_{i_1} T_{i_2} \cdots T_{i_l}$.  It
is a consequence of the relations (\ref{heckegen-2})-(\ref{heckegen-3})
that $T_{\sigma}$ depends only on $\sigma$, and not on the factorization
into adjacent transpositions.

Let $z_1^{d_1} v_{m_1} \otimes z_2^{d_2} v_{m_2}
\otimes \cdots \ $ be a semi-infinite tensor.
Define its {\em $q$-antisymmetrization} to be
\begin{equation}
z^{d_1} v_{m_1} \wedge_q z^{d_2} v_{m_2} \wedge_q \cdots = \sum_{\sigma
\in S_{\infty}} (z_1^{d_1} v_{m_1} \otimes z_2^{d_2} v_{m_2} \otimes
\cdots ) \cdot T_{\sigma}. \label{hecke-antisym}
\end{equation}
Here the action of $T_{\sigma}$ is given by equation (\ref{hecke-action}).
Again, subscripts on the variables have been dropped in infinite wedge
notation.  By the first remark above, (\ref{hecke-antisym}) really is an
antisymmetrization, rather than a symmetrization.

\bigskip
\noindent
{\bf Conjecture} \hspace{1mm}
The sum given by (\ref{hecke-antisym}) converges in the power series
topology.  (I.e., the coefficient of each particular tensor $z_1^{j_1}
v_{k_1} \otimes z_2^{j_2} v_{k_2} \otimes \cdots $ is a well-defined
power series in $q$.)

\medskip
The formal action of $K_i$, $E_i$, $F_i
\in U_q(\widehat{\frak{sl}}_n)$ on semi-infinite
tensors gives rise to a genuine action of $U_q(\widehat{\frak{sl}}_n)$
on the vector space spanned the $q$-wedges.  As in Section 2, this vector
space is a level 1 module.  Since each $T_{\sigma}$ commutes with the
action of $U_q(\widehat{\frak{sl}}_n)$, when $X \in U_q(\widehat{
\frak{sl}}_n)$ acts on the $q$-antisymmetrization of $z_1^{d_1} v_{m_1}
\otimes z_2^{d_2} v_{m_2} \otimes \cdots \ $, the result is the
$q$-antisymmetrization of $\Delta^{\infty}(X) \cdot (z_1^{d_1} v_{m_1}
\otimes z_2^{d_2} v_{m_2} \otimes \cdots )$.
\begin{prop}
The $q$-wedge
$$v_{\Lambda_i} = zv_i \wedge_q zv_{i-1} \wedge_q \cdots \wedge_q
zv_1 \wedge_q z^2v_n \wedge_q \cdots \wedge_q z^2v_1 \wedge_q \cdots$$
is a highest weight vector of weight $\Lambda_i$.
\end{prop}
{\em Proof.} \hspace{2mm}  By the above, it is enough to check that
if $z_i^{d_i} v_{m_i} \otimes z_{i+1}^{d_i} v_{m_i}$ appears somewhere
in a tensor, then the $q$-antisymmetrization of that tensor
is zero.  Grouping the
basis elements of $H_{\infty}(q^2)$ in pairs $\{T_{\sigma}, T_i
T_{\sigma} \}$ corresponding to cosets of the
subgroup $H=\{ \mbox{id}, \sigma_i \} \subset S_{\infty}$ reduces the
problem to showing that such a vector is killed by $1 + T_{\sigma_i}$.
This is an immediate consequence of the first remark following
equation (\ref{hecke-action}).

\medskip
The highest weight vector $v_{\Lambda_i}$ generates a highest weight
module $V_{\Lambda_i}$ of weight $\Lambda_i$
inside the space $F_{(i)}$ spanned by
the $q$-wedges that are the same as $v_{\Lambda_i}$ after
finitely many terms.  $V_{\Lambda_i}$ is a $q$-deformation of the
corresponding $\widehat{\frak{sl}}_n$-module constructed in Section
\ref{affine-class}, and is therefore irreducible.

\bigskip
\noindent
{\bf Remarks}

\smallskip
\noindent
1.  The $q$-antisymmetrization for $U_q(\frak{sl}_{\infty})$
given by (\ref{q-antisym}) is essentially a special case of the one
for $U_q(\widehat{\frak{sl}}_n)$ given by (\ref{hecke-antisym}), since
it follows from (\ref{hecke-action}) that
$$ (v_{k_1} \otimes v_{k_2} \otimes \cdots ) \cdot T_i =
-q v_{k_1} \otimes \cdots \otimes v_{k_{i+1}} \otimes v_{k_i} \otimes
v_{k_{i+2}} \otimes \cdots $$
as long as $k_i > k_{i+1}$.

\smallskip
\noindent
2.  The action of $U_q(\widehat{\frak{sl}}_n)$ on $F_{(i)}$ was
originally constructed by Hayashi in \cite{Hayashi}.

\medskip
The following lemma relating the $q$-antisymmetrizations of $z_1^{j_1}
v_{m_1} \otimes z_2^{j_2} v_{m_2} \otimes \cdots \ $ and of $(z_1^{j_1}
v_{m_1} \otimes z_2^{j_2} v_{m_2} \otimes \cdots ) \cdot T_i$ will be
useful in the next section:

\medskip
\noindent
{\bf Lemma}  \hspace{1mm} {\em For any $T_i$,}
$$\sum_{\sigma \in S_{\infty}} (z_1^{j_1} v_{m_1} \otimes z_2^{j_2} v_{m_2}
\otimes \cdots ) T_i \cdot T_{\sigma} = q^2 \sum_{\sigma \in S_{\infty}}
(z_1^{j_1} v_{m_1} \otimes z_2^{j_2} v_{m_2} \otimes \cdots ) \cdot
T_{\sigma}. $$
\medskip
\noindent
{\em Proof of Lemma.} \hspace{2mm}
As before, group the basis elements of $H_{\infty}(q^2)$
in pairs $\{T_{\sigma} , T_i T_{\sigma} \}$ corresponding to left
cosets of $\{ \mbox{id}, \sigma_i \}$ in $S_{\infty}$.  Using
the first Hecke relation (\ref{heckegen-1}), compute
\begin{eqnarray*}
T_i (T_{\sigma} + T_i T_{\sigma} ) &=& T_i T_{\sigma} + ((q^2 -1) T_i
+ q^2 ) T_{\sigma} \\
&=& q^2 (T_i T_{\sigma} + T_{\sigma}).
\end{eqnarray*}
The lemma follows.

\subsection{Vertex operators}
Splitting off the first component of every tensor defines intertwiners
$$\widetilde{\Phi}_{(i)} : F_{(i)} \to V(z) \
\widehat{\otimes} \ F_{(i-1)}$$
as in Section \ref{vertex}.  Since $F_{(i)}$ is a $q$-deformation of the
completely reducible $\widehat{\frak{sl}}_n$-module $F_{(i)}$ from Section
2, it is also completely reducible.  This makes it possible to define
intertwiners $\Phi_{(i)} : V_{\Lambda_i} \to V(z) \ \widehat{\otimes}
\ V_{\Lambda_{i-1}}$ by composing with the projection $F_{(i-1)} \to
V_{\Lambda_{i-1}}$.
These intertwiners satisfy the obvious
$q$-analog of Proposition \ref{vertex-comp}.  Up to normalization,
they are the vertex operators studied in \cite{qKZ} and \cite{vertex}
by means of the quantum KZ equation.  Their iterations can be computed
exactly as in the classical case.

As an example, take the case of $U_q(\widehat{\frak{sl}}_2)$ and
consider the operator $\widetilde{\Phi}_{(0)}: V_{\Lambda_0} \subset
F_{(0)} \to V(z) \ \widehat{\otimes} \ F_{(1)}$.
\begin{prop}
The image of $v_{\Lambda_0}$ under $\widetilde{\Phi}_{(0)}$ is given by
\begin{eqnarray}
\widetilde{\Phi}_{(0)} (v_{\Lambda_0}) &=&
\sum_{j=1}^{\infty} q^{3(j-1)} z^j v_2 \otimes (zv_2 \wedge_q zv_1 \wedge_q
\cdots \wedge_q z^{j-1}v_1 \wedge_q z^j v_1 \wedge_q z^{j+1} v_2 \wedge_q
\cdots)  \label{act-high-wt} \\
& & - \sum_{j=1}^{\infty} q^{3(j-1)+1} z^j v_1 \otimes (zv_2 \wedge_q
zv_1 \wedge_q \cdots \wedge_q z^j v_2 \wedge_q z^{j+1} v_2 \wedge_q
z^{j+1} v_1 \wedge_q \cdots). \nonumber
\end{eqnarray}
\end{prop}
{\em Proof.}  \hspace{2mm}
The assertion states that the $q$-antisymmetrization
of $z v_2 \otimes z v_1 \otimes z^2 v_2 \otimes \cdots \ $ is given
by equation (\ref{act-high-wt}).  (For simplicity of notation, the
subscripts on the variables have been left off, but again it should be
understood that all $z$'s have subscripts on them recording the factor
in which they appear.)  The idea of the proof is to compute the
antisymmetrization using
(\ref{hecke-action}).  Each term in (\ref{act-high-wt}) corresponds to
a collection of tensors with
a particular basis element $z^j v_k$ appearing in the first component.
So imagine performing the antisymmetrization by first moving into the
first component whichever $z^j v_k$ will go there: choose a particular
$z^j v_k$ (say the one in the $r$-th component) and apply $T_{r-1},
T_{r-2}, \ldots , T_1$ to move it all the way over to the left.  The
effect of one step in this sequence
can be computed using the following formulas
for the action of $H_2(q^2)$ (with generator $T$) on $V(z) \otimes
V(z)$ (again, it should be understood that the $z$'s in the first and
second components represent different $z_i$'s):
\begin{eqnarray}
(z^j v_i \otimes z^k v_i) \cdot T &=& - q^2 z^k v_i \otimes z^j v_i \ - \
(q^2-1) (z^{k-1} v_i \otimes z^{j+1} v_i \ + \ z^{k-2} v_i \otimes z^{j+2}
v_i \label{nub-1} \\
& & \hspace{2in} + \cdots + \ z^{j+1} v_i \otimes z^{k-1} v_i) \nonumber \\
(z^j v_1 \otimes z^k v_2) \cdot T &=& - q z^k v_2
\otimes z^j v_1 \ - \ (q^2-1)
(z^{k-1} v_1 \otimes z^{j+1} v_2 \ +
\ z^{k-2} v_1 \otimes z^{j+2} v_2
\label{nub-2} \\
& & \hspace{2in} + \cdots + \ z^{j+1} v_1
\otimes z^{k-1} v_2) \nonumber \\
(z^j v_2 \otimes z^k v_1) \cdot T
&=& - q z^k v_1 \otimes z^j v_2 \ - \ (q^2-1)
(z^k v_2 \otimes z^j v_1 \ +
\ z^{k-1}v_2 \otimes z^{j+1} v_1 \label{nub-3} \\
& & \hspace{2in} + \cdots + \ z^{j+1} v_2 \otimes z^{k-1} v_1 )
\nonumber
\end{eqnarray}
These formulas follow directly from (\ref{hecke-action}). Here $j<k$
except in the last equation, where $j \leq k$.  The upshot is that after
applying a sequence $T_{r-1}, T_{r-2}, \ldots , T_1$, the result will be
\begin{eqnarray}
\lefteqn{
\sum_{j=1}^{\infty} q^{3(j-1)} z^j v_2 \otimes (zv_2 \otimes zv_1 \otimes
\cdots \otimes z^{j-1}v_1 \otimes z^j v_1 \otimes z^{j+1} v_2 \otimes
\cdots) } \label{move-first} \\
& & \hspace{4mm}
- \sum_{j=1}^{\infty} q^{3(j-1)+1} z^j v_1 \otimes (zv_2 \otimes
zv_1 \otimes \cdots \otimes z^j v_2 \otimes z^{j+1} v_2 \otimes
z^{j+1} v_1 \otimes \cdots) + \ \mbox{other terms,} \nonumber
\end{eqnarray}
where the other terms, which come from the terms following the first one
in (\ref{nub-1})-(\ref{nub-3}), all have to the right of the first component
one of the following four sequences:
\begin{eqnarray}
z^k v_1 \otimes z^j v_1 \otimes z^{j+1} v_2 \otimes z^{j+1} v_1 \otimes
\cdots \otimes z^{k-1} v_1 \otimes z^k v_2 \otimes z^k v_1
\label{kill-1} \\
z^k v_1 \otimes z^j v_2 \otimes z^j v_1 \otimes z^{j+1} v_2 \otimes
\cdots \otimes z^{k-1} v_1 \otimes z^k v_2 \otimes z^k v_1
\label{kill-2} \\
z^k v_2 \otimes z^j v_1 \otimes z^{j+1} v_2 \otimes z^{j+1} v_1 \otimes
\cdots \otimes z^{k-1} v_2 \otimes z^{k-1} v_1 \otimes z^k v_2
\label{kill-3} \\
z^k v_2 \otimes z^j v_2 \otimes z^j v_1 \otimes z^{j+1} v_2 \otimes
\cdots \otimes z^{k-1} v_2 \otimes z^{k-1} v_1 \otimes z^k v_2
\label{kill-4}
\end{eqnarray}
Here $j \leq k$, and in the third sequence, $j<k$.

At this point, the first component is left alone, and the other components
are completely $q$-antisymmetrized.  Antisymmetrizing the terms listed
explicitly in (\ref{move-first}) gives the answer.  The other terms
do not figure in the answer because strings of the form
(\ref{kill-1})-(\ref{kill-4}) antisymmetrize to $0$.
This can be seen by induction on the length of
the string, using the lemma in the preceding section.  If the rightmost
$z^k v_1$ (or $z^k v_2$, as appropriate), appears in the $r$-th
component, apply $T_{r-1}, T_{r-2}, \ldots $ to move it over next
to the $z^k v_1$ or $z^k v_2$ on the left.  A term containing a
$z^k v_1 \otimes z^k v_1$ or a $z^k v_2 \otimes z^k v_2$ antisymmetrizes
to $0$, and by the lemma, it is enough to show that the extra terms
created by this procedure also antisymmetrize to $0$.
By (\ref{nub-1})-(\ref{nub-3}), all these extra terms contain sequences
of the form (\ref{kill-1})-(\ref{kill-4}) of shorter length, which
antisymmetrize to $0$ by induction.

\bigskip
\noindent
{\bf Remarks}

\smallskip
\noindent
1.  The formula for $\widetilde{\Phi}_{(0)} (v_{\Lambda_0})$ given by
equation (\ref{act-high-wt}) looks (perhaps with all $z_i$ set equal to
$1$) exactly like formula (A3.1) in \cite{vertex} giving the expansion of
$\Phi_{(0)} (v_{\Lambda_0})$ in terms of the upper global base for
$V_{\Lambda_1}$.  This suggests that $q$-wedges may be related in a
reasonable way to the upper global base.  (Notice that the
coproduct used in \cite{vertex} differs by a permutation from the one
given here, and the semi-infinite tensor product considered there is
$\cdots V(z_3) \otimes V(z_2) \otimes V(z_1)$.  Therefore, the tensors
appearing in \cite{vertex} should be reversed to match the conventions
of this paper.)

\smallskip
\noindent
2.  As in the classical case, solutions to approriate $q$-KZ equations
(for example, those
listed on p. 116 of \cite{vertex}) should be understood (up to scalars)
as finite $q$-antisymmetrizations.  The functions $w_n(z)$ (see pp.
121-122 of \cite{vertex}) are (the reverses of) the $q$-wedges
$v_2 \wedge_q v_1 \wedge_q \cdots \wedge_q z^{m-1} v_2 \wedge_q
z^{m-1} v_1$ if $n=2m$, and $v_1 \wedge_q z v_2 \wedge_q \cdots
\wedge_q z^{m-1} v_2 \wedge_q z^{m-1} v_1$ if $n=2m-1$.
Also as in the classical case,
the infinite $q$-antisymmetrization $v_{\Lambda_0}$
should be nothing more than a limit of these solutions as the number of
variables goes to infinity.  The details will be worked out in a separate
publication.

\smallskip
\noindent
3.  In \cite{vertex}, a tensor product $\cdots \otimes V(z_i) \otimes
V(z_{i+1}) \otimes \cdots$ which is infinite in both directions is
considered.  The problem is to find the ground state vector of the XXZ
Hamiltonian, which spans the trivial component of $V(\Lambda_0)^* \otimes
V(\Lambda_0) \subset \cdots \otimes V(z_i) \otimes
V(z_{i+1}) \otimes \cdots$.  From the point of view outlined here,
this vector should be a $q$-antisymmetrization
which is infinite in both directions.

\vspace{.6in}
\noindent
I am grateful to Nicolai Reshetikhin for introducing me to this
subject and to this problem, and for many stimulating conversations.
I am also grateful to Tetsuji Miwa for reading an early draft of
this paper and pointing out inaccuracies.

\vspace{5mm}
\noindent
{\sc Department of Mathematics \\
\noindent
University of California at Berkeley \\
\noindent
Berkeley, CA 94720} \\
\noindent
{\tt stern@math.berkeley.edu}
\end{document}